\documentclass[12pt]{article}

\usepackage{amsfonts}
\usepackage{amsmath,amssymb}
\usepackage{amscd}
\usepackage{amsthm}
\newtheorem{df}{Definition}[section]

\newtheorem{lem}{Lemma}[section]
\newtheorem{cly}{Corollary}[section]
\newtheorem{prop}{Proposition}[section]
\newtheorem{rem}{Remark}[section]
\newtheorem{ex}{Example}[section]

\title{On Hom type algebras}
\author {Ya\"el Fr\'egier; Aron Gohr\\{\small Mathematics Research Unit}\\
{\small  162A, avenue de la faiencerie}\\
{\small  L-1511, Luxembourg}\\
 {\small  Grand-Duchy of Luxembourg.}\\
}
\begin {document}
\date{}
\maketitle


\begin{abstract}

Hom-algebras are generalizations of algebras obtained using a
twisting by a linear map. But there is a priori a freedom on where
to twist. We enumerate here all the possible choices in the Lie and
associative categories and study the relations between the obtained
algebras. The associative case is richer since it admits the notion
of unit element. We use this fact to find sufficient conditions for
hom-associative algebras to be associative and classify the
implications between the hom-associative types of unital algebras.
\end{abstract}



\section* { Introduction.}

The present paper investigates variations on the theme of \emph{Hom-algebras}, a topic which has recently received much attention from various researchers \cite{CaeGoy}, \cite{Goh09}, \cite{Yau:HOMology}, \cite{Yau:Hom-quaI}, \cite{Yau:Hom-quaII}. Generally speaking, the notion of Hom-algebra over a certain operad is obtained by twisting in a strategic way the identities for the algebra multiplication implied by the operad in question. For instance, an algebra $(V, \star)$ together with a linear self-map $\alpha:V \rightarrow V$ is called \emph{Hom-associative} if it satisfies the identity
\begin{displaymath}
\mu \circ (\alpha \otimes \mu) = \mu \circ (\mu \otimes \alpha)
\end{displaymath}
which is obtained by replacing $id$ with $\alpha$ in the ordinary associativity condition
\begin{displaymath}
\mu \circ (id \otimes \mu) = \mu \circ (\mu \otimes id).
\end{displaymath}
The study of Hom-associative algebras originates with work by Hartwig, Larson and Silvestrov in the Lie case \cite{HLS}, where a notion of Hom-Lie algebra was introduced in the context of studying deformations of Witt and Virasoro algebras. Later, it was extended to the associative case by Makhlouf and Silverstrov in \cite{MSI}. A number of classical constructions have been found to have a Hom-counterpart, see e.g. \cite{AMS}, \cite{LS1}, \cite{MSII},
\cite{MSIII}, \cite{MSIV}, \cite{Yau:EnvLieAlg}, \cite{Yau:HOMology}
and \cite{Yau:Hom-bi}.

When studying the structure theory of hom-associative algebras, one naturally encounters related algebras which, while not obeying the traditional hom-associative identities, satisfy similar identities \cite{FGII}. For instance, in the proof of the main result of \cite{FGII}, there naturally appear algebraic structures $(V, \star, \alpha)$ which satisfy the conditions
\begin{displaymath}
(\alpha(x) \star y) \star z = x \star (y \star \alpha(z))
\end{displaymath}
or
\begin{displaymath}
\alpha(x \star y) \star z = x \star \alpha(y \star z)
\end{displaymath}
for all $x,y,z \in V$. This observation suggests that in the context of studying Hom-algebraic structures, there is a natural interest in exploring alternative possibilities on ``how to twist'' the identities of a given classical algebraic category to obtain a Hom-counterpart.\\
In this paper, we start a systematic exploration of other possibilities to define Hom-type algebras. We will, for the purposes of the present paper, restrict the scope of this first investigation in several ways. First, we will only consider Hom-twisted versions of the associative and Lie categories. Second, we will only consider symmetric twisted identities and limit the ``degree'' in terms of occurrences of the twisting map $\alpha$ in the defining identity of the twisted categories.\\
The paper is partitioned in two main sections. In the first section we introduce, following the ideas outlined above, new types of Hom-Lie algebras and study their properties in special cases. We give several examples of these new types of Hom-Lie algebras and study their relations among each other and to ordinary Lie algebras. In the second section, we introduce in analogy to this work a similar system of new types of Hom-associative algebras. We point out that in the case of unital algebras, these types of Hom-associativity conditions can be partially ordered by restrictiveness, with the traditional hom-associativity condition ending up on top, i.e. as most restrictive. Finally, we introduce \emph{Hom-monoids} to obtain an easy way to construct counterexamples to possible relations between types of Hom-algebras which do not hold. These counterexamples prove that our partial ordering of hom-type algebras cannot be improved upon.\\
We end the introduction by fixing some conventions and notations. In this article, $k$ will by default be a commutative ring, $K$ a field. Modules and algebras will by default be understood to be over an arbitrary commutative ring. If $\alpha:G \rightarrow H$ is a homomorphism of groups (rings, modules, etc.) we will denote by $Ke(\alpha)$ its kernel and by $Im(\alpha)$ its image. $V$ will by default be a $k$-module.

\section { Hom-Lie algebras.}

In this part we define types of hom-Lie algebras and give some
relations between them.

\subsection { Definitions.} \label{liedef}

We start by recalling the original definition  following
\cite{HLS}.

\begin{df}
\label{defhomlie}
A Hom-Lie algebra is a triple $(V, [\cdot, \cdot], \alpha)$
consisting of a module $V$ over a commutative ring $k$, a bilinear map $[\cdot, \cdot]$:
$V\times V\longrightarrow V$ and a linear space homomorphism $\alpha
: V\longrightarrow V  $ satisfying
\begin{eqnarray}
[x,x] = 0\\
\circlearrowleft_{x,y,z} [\alpha (x),[y,z]]=0 \label{I_1}
\end{eqnarray}
for all $x,y,z$ in V, where $\circlearrowleft_{x,y,z}$ denotes
summation over the cyclic permutations on $x,y,z$. Explicitly, this means
\begin{displaymath}
\circlearrowleft_{x,y,z} [\alpha(x),[y,z]] := [\alpha(x),[y,z]] + [\alpha(y), [z,x]] + [\alpha(z), [x,y]].
\end{displaymath}
\end{df}
Note that, if $(G, +, 0)$ is a group with identity $0$ and $\varphi: G \times G \rightarrow G$ is a bi-additive map with $\varphi(x,x) = 0$, we automatically have $\varphi(x,y) + \varphi(y,x) = 0$. On the other hand, if in $G$ the equation $x + x = 0$ always implies $x = 0$, as is the case for instance if $G$ is the additive group of a field of characteristic $\neq 2$, then the condition $\varphi(x,y) + \varphi(y,x) = 0$ implies also $\varphi(x,x) = 0$ for all $x \in G$. Therefore, condition $(1)$ in our definitions corresponds to the usual condition of skew-symmetry of a Lie bracket.\\

Now if we look at (\ref{I_1}), it is natural to ask why we chose to twist by
$\alpha$ in the first argument, and not in the second or third? This question is the first motivation for us to suggest the introduction of two new types, $I_2$ and $I_3$, of Hom-Lie algebras:

\begin{df}
 A Hom-Lie algebra of type $I_2$ is defined by
replacing, in Definition \ref{defhomlie}, equation (\ref{I_1}) by
\begin{eqnarray}
\circlearrowleft_{x,y,z} [x,[\alpha (y),z]]=0 \label{I_2}.
\end{eqnarray}
 If one uses
\begin{eqnarray}
\circlearrowleft_{x,y,z} [x,[y,\alpha (z)]]=0 \label{I_3}
\end{eqnarray}
instead, one gets the definition of a Hom-Lie algebra of type $I_3$.
\end{df}

  A Hom-Lie algebra in the usual sense should be referred to as "Hom-Lie
algebra of type $I_1$", but we will, most of the time, simply use
the term "Hom-Lie algebra", for coherence with the usage in the
literature.

\begin{rem}\label{rem1} Of course type $I_2$ and $I_3$ are the same by skew symmetry of the
bracket. Nevertheless we introduce these two types for pedagogical
reasons, since they will appear again in the associative category.
\end{rem}

Now we remark that $\alpha$ has still two more choices for its
dinner: it could be applied to the results of the first or of the
second bracket, and give two other types.

\begin{df}
 If one replaces in Definition \ref{defhomlie}, equation (\ref{I_1})
by
\begin{eqnarray}
\circlearrowleft_{x,y,z} [x,\alpha ([y,z])]=0 \label{II}
\end{eqnarray}

one gets the definition of a Hom-Lie algebra of type $II$.\\
 If one uses
\begin{eqnarray}
\circlearrowleft_{x,y,z} \alpha([x,[y,z]])=0 \label{III},
\end{eqnarray}
instead, one gets the definition of a Hom-Lie algebra of type $III$.
\end{df}

A trivial example of an algebra of type $III$ is given by
considering an arbitrary Lie algebra structure on $V$ (a Hom-Lie
algebra where the twisting is the identity on V), together with an
arbitrary linear $\alpha$. In the other direction, we see that if $\alpha$ is any \emph{injective} linear map and $(V, \alpha)$ is a Hom-Lie algebra of type $III$, then $V$ must be Lie.\\
 Another example, less trivial, is obtained considering the notion of descending central series $V^n$,
borrowed from Lie theory: $V^0:=V, V^1:=[V,V], V^2:=[V,V^1],
\dots,V^n:=[V,V^{n-1}]$. Considering an arbitrary $\alpha$ whose
kernel contains $V^2$ gives the second example. One can also obtain examples where $\alpha$ does not vanish on $V^2$. The following example is an extreme case of this insofar as the kernel of $\alpha$ is one-dimensional, i.e. of lowest possible dimension:
\begin{ex}
Let $K$ be a field and let $V := K^3$. We define a bilinear map $[\cdot, \cdot]: V \times V \rightarrow V$ by
\begin{displaymath}
\left[ \begin{pmatrix} \lambda_1 \\ \lambda_2 \\ \lambda_3\\ \end{pmatrix}, \begin{pmatrix} \mu_1\\ \mu_2 \\ \mu_3\\ \end{pmatrix} \right] := \begin{pmatrix} \lambda_1 \mu_3 - \lambda_3 \mu_1 \\ 0 \\ \lambda_2 \mu_3 - \lambda_3 \mu_2 \end{pmatrix}.
\end{displaymath}
It is clear that this map satisfies $[v,v] = 0$ for all $v \in V$. Also, our bracket does not induce a structure of Lie algebra on $V$, since e.g. with $e_1, e_2, e_3$ the canonical basis vectors of $V$ we have
\begin{displaymath}
\circlearrowleft_{e_1, e_2, e_3} [e_1, [e_2, e_3]] = e_1 \neq 0.
\end{displaymath}
Finally, with $\alpha: V \rightarrow V$ defined through
\begin{displaymath}
\alpha \begin{pmatrix} \lambda_1 \\ \lambda_2 \\ \lambda_3 \end{pmatrix} := \begin{pmatrix} \lambda_2 \\ \lambda_3 \\ 0 \end{pmatrix}
\end{displaymath}
we see that $(V, \alpha, [\cdot, \cdot])$ is Hom-$III$-Lie by a straightforward calculation. But $V^2$ has as basis the set $\{e_1, e_3\}$, so $\alpha$ does not vanish on $V^2$.
\end{ex}

Let us now do a little science fiction and imagine that $\alpha$ is
in fact a morphism of algebra, i.e satisfies
$[\alpha(x),\alpha(x)]=\alpha([x,y]) \ \forall x,y \in V$. One could
then, using this equality twice,  rewrite equation (\ref{III}) in
two ways ((\ref{III'}) and (\ref{III''})). Hence it is natural to
also consider:

\begin{df}
 Hom-Lie algebras of types  $III'$ and  $III''$ are defined by respectively replacing
in Definition \ref{defhomlie}, equation (\ref{I_1}) by
\begin{eqnarray}
 \circlearrowleft_{x,y,z} [\alpha(x),\alpha ([y,z])]=0,
\label{III'}\\
 \circlearrowleft_{x,y,z} [\alpha(x),[\alpha
(y),\alpha(z)]]=0 \label{III''}.
\end{eqnarray}
\end{df}

One can remark that a Hom-Lie algebra of type  $III''$ is nothing
else than a Lie algebra structure on the image of $\alpha$.\\

Similarly, (\ref{II}) leads to the equation (\ref{II_1}) which is
quadratic in $\alpha$. But once we have opened the Pandora's
box, we are forced to also consider the other quadratic expressions
in $\alpha$, (\ref{II_2}) and (\ref{II_3}):

\begin{df}
 Hom-Lie algebras of types $II_1$, $II_2$ and  $II_3$ are defined by respectively replacing
in Definition \ref{defhomlie}, equation (\ref{I_1}) by
\begin{eqnarray}
\circlearrowleft_{x,y,z} [x,[\alpha (y),\alpha(z)]]=0 \label{II_1},\\
 \circlearrowleft_{x,y,z} [\alpha(x),[y,\alpha (z)]]=0\label{II_2},\\
\circlearrowleft_{x,y,z} [\alpha(x),[\alpha (y),z]]=0 \label{II_3}.
\end{eqnarray}
\end{df}

\begin{rem}\label{rem2} It is easy to see that Remark \ref{rem1} applies mutatis mutandis if one replaces $I_2$ and $I_3$ by $II_2$ and $II_3$.
 \end{rem}

\subsubsection*{On the notations:}

 We need a notation to distinguish all these types of
Hom-Lie algebras.

 "$Hom^{\text{type}}-Lie$ algebra" seems appropriate.
For example "$Hom^{I_2}-Lie$ algebra" will stand for "Hom-Lie
algebra of type $I_2$". By $Hom^\star-Lie$ we will mean a
Hom-algebra of simultaneously all types.

 We have chosen to divide these
types in three classes $I$, $II$ and $III$ accordingly to the degree
in $\alpha$, i.e. the number of occurrences of $\alpha$ in the
defining equations. We consider the "virtual" degree. In
$\alpha([x,y])$ for example, $\alpha$ is of virtual degree two even
if it appears only once. This is because if $\alpha$ is a morphism
for the bracket, one has $\alpha([x,y])=[\alpha(x),\alpha(x)]$ which
is really of degree two in $\alpha$. We hope that this choice will
help the reader to memorize easily the subdivision in classes.

 We used the "prime" notation to remember that $III'$ is derived from
 $III$, and that $III''$ is derived from $III'$. Accordingly, $II_1$, $II_2$ and $II_3$ should have been denoted by  $II_1'$, $II_2'$ and
 $II_3'$, but we decided to omit the upper-script prime, since the
 lower script enables already to distinguish these types.

We have chosen the ordering in the classes II and I in a way that
they coincide under the symmetry $S$ which consists of interchanging
the role of $\alpha$ and $id$ (the identity of $V$), namely for example
$S([Id(x),[\alpha (y),\alpha(z)]]):=[\alpha (x),[Id(y),Id(z)]]$.

Finally we introduce the Jacobiator associated to each of these
structures as the left hand side of the defining equation of the $Hom^{type}-Lie$ algebra under consideration. One denotes it by
$J^{type}_\alpha$. As an example,
$$J^{II_1}_\alpha:=\circlearrowleft_{x,y,z} [x,[\alpha
(y),\alpha(z)]].$$ In particular $S(J^{I_i}_\alpha)=J^{II_i}_\alpha$
for $1\leq i\leq 3$.

We conclude this part by the following table which summarizes the
list of types of Hom-Lie algebras:
\\

\begin{tabular}{|c|c|}
\hline $I_1$  & $\circlearrowleft_{x,y,z}[\alpha (x),[y, z]]=0$
\\
\hline $I_2$  & $\circlearrowleft_{x,y,z}[x,[\alpha (y), z]]=0$\\
\hline $I_3$  & $\circlearrowleft_{x,y,z}[x,[y, \alpha (z)]]=0$\\
 \hline
\end{tabular}

\begin{tabular}{|c|c|}
\hline $II_1$  & $\circlearrowleft_{x,y,z}[x,[\alpha (y), \alpha
(z)]]=0$
\\
\hline $II_2$  & $\circlearrowleft_{x,y,z}[\alpha (x),[y, \alpha (z)]]=0$\\
\hline $II_3$  & $\circlearrowleft_{x,y,z}[\alpha (x),[\alpha (y),z]]=0$\\
 \hline
\end{tabular}

\begin{tabular}{|c|c|}
\hline $III$  & $\circlearrowleft_{x,y,z}\alpha ([x,[y, z]])=0$\\
\hline $III'$  & $\circlearrowleft_{x,y,z}[\alpha (x),\alpha([y,z])]=0$\\
 \hline
 \hline $III''$  & $\circlearrowleft_{x,y,z}[\alpha(x),[\alpha (y),
\alpha (z)]]=0$
\\
\hline
\end{tabular}

\subsection { Relations among these types.} \label{lierel}

Now that we have these new types of Hom-Lie algebras, it is natural
to seek for relations among them.

\begin{prop}\label{prop1}
Let us suppose that $(V, [ \cdot , \cdot ])$ is a Lie algebra and consider
 $(V, [ \cdot , \cdot ], \alpha)$,  
\begin{enumerate}
\item if it is a hom-Lie algebra of type $I_2$, it is necessarily also of type $I_1,$
\item if it is a hom-Lie algebra of type $II_2$, it is necessarily of type $II_1.$
\end{enumerate}
\end{prop}
In particular
\begin{cly}
 If moreover $\alpha$ is a morphism of $[ \cdot , \cdot ]$, being of type
$II_2$ implies to be of type  $II_1$ which in turn is equivalent to be of type $II$.
\end{cly}

\begin{proof} We start by proving the assertion that a hom-Lie algebra of type $I_2$ is necessarily  of type $I_1$. Let us first establish the following property:
\begin{equation}\label{property}
J_{\alpha}^{I_1}=-J_{\alpha}^{I_2}-J_{\alpha}^{I_3}.
 \end{equation}
Indeed, since the bracket satisfies the Jacobi identity, one has:
$$[\alpha (x),[y,z]]=
-[z,[\alpha
(x),y]]-[y,[z,\alpha (x)]].$$

Summing this relation over cyclic permutations leads to the desired property.

Now let us suppose that we have a $Hom^{I_2}-Lie$ algebra, i.e.
$J_{\alpha}^{I_2}=0$. Remark \ref{rem2} implies that one also has
$J_{\alpha}^{I_3}=0$. The property (\ref{property}) enables to
conclude.

The proof that a hom-Lie algebra of type $II_2$ is necessarily of type $II_1$ is almost
the same, one just needs to read the preceding proof after
having applied the symmetry $S$.
\end{proof}

The reverse implication in the preceding proposition would need, to
hold, that $J_{\alpha}^{I_2}+J_{\alpha}^{I_3}=0 \Rightarrow
J_{\alpha}^{I_2}=J_{\alpha}^{I_3}=0$. It is natural to ask for a
necessary and sufficient condition for this last implication. We do
not know the answer to this question but we can remark the ''self-adjointness'' condition $[\alpha
(x),y]=[x,\alpha (y)] \ \forall x,y\in V$ on $\alpha$ is sufficient if the underlying $k$-module $V$ is $2$-torsion-free. But this condition of self-adjointness is fairly strong, implying for instance that $[\alpha(x), x] = 0$ for all $x \in V$ if $V$ is a vector space over a field of characteristic $\neq 2$, and is therefore in most cases unsatisfied. The following example shows that $J_{\alpha}^{I_2} + J_{\alpha}^{I_3} = 0$ does not indeed imply $J_{\alpha}^{I_2} = 0$, as would be expected:
\begin{ex}
Let $K$ be a field with $char(K) \neq 2$ and let $V := K^2$. We define on $V$ the bracket
\begin{displaymath}
\left[ \begin{pmatrix}\lambda_1 \\ \lambda_2 \end{pmatrix}, \begin{pmatrix} \mu_1 \\ \mu_2 \end{pmatrix} \right] := \begin{pmatrix} 0 \\ \lambda_1 \mu_2 - \lambda_2 \mu_1 \end{pmatrix}.
\end{displaymath}
It is clear that this is skew-symmetric, and direct calculation verifies the Lie identity. Set now
\begin{displaymath}
\alpha \begin{pmatrix} \lambda_1 \\ \lambda_2 \end{pmatrix} := \begin{pmatrix} \lambda_1 + \lambda_2 \\ \lambda_2 \end{pmatrix}.
\end{displaymath}
Then we see
\begin{displaymath}
\left[ \alpha \begin{pmatrix} a_1 \\ a_2 \end{pmatrix}, \left[ \begin{pmatrix} b_1 \\ b_2 \end{pmatrix}, \begin{pmatrix} c_1 \\ c_2 \end{pmatrix} \right] \right] = \begin{pmatrix} 0 \\ a_1 b_1 c_2 + a_2 b_1 c_2 - a_1 b_2 c_1 - a_2 c_1 b_2 \end{pmatrix}
\end{displaymath}
and summing up cyclic permutations of the term in the second component this implies that condition $J_{\alpha}^{I_1}$ is satisfied. On the other hand we have $[e_1, [\alpha(e_2), e_2]] \neq 0$, so $(V, [\cdot, \cdot])$ is not of type $I_2$.
\end{ex}

Let again $V$ be an arbitrary $k$-module, suppose that $[\cdot, \cdot]: V \times V \rightarrow V$ is a bilinear map and suppose that $\alpha:V \rightarrow V$ is linear. We consider now on V the bracket $[\cdot, \cdot]_\alpha$ defined by
$[a,b]_\alpha:= [a,b]+[\alpha (a),b]+ [a,\alpha (b)]$. We have the
following:

\begin{prop}\label{prop2}
Let us suppose that $(V, [ \cdot , \cdot ],\alpha)$ is a
$Hom^{\star}-Lie$, then $[ \cdot , \cdot ]$ is a Lie algebra if and only if $[\cdot,
\cdot]_\alpha$ is a Lie algebra.
\end{prop}

\begin{proof}
Let us introduce the notation $\mathfrak{J}_\alpha$, (resp
$\mathfrak{J}$) for the jacobiator of $[\cdot, \cdot]_\alpha$ (resp
$[\cdot, \cdot]$). A direct computation gives that
\begin{equation}\label{lie_alpha}
\mathfrak{J}_\alpha=\mathfrak{J}+
J^{I_1}_\alpha+J^{I_2}_\alpha+J^{I_3}_\alpha+J^{II}_\alpha+J^{II_2}_\alpha+J^{II_3}_\alpha.
\end{equation}
Being a $Hom^\star$-Lie algebra implies that this last equality
reduces to
$$\mathfrak{J}_\alpha=\mathfrak{J},$$
hence the conclusion.
\end{proof}

\begin{rem}
We did not use the hypothesis of being a $Hom-Lie$ algebra of types
$III$, $III'$, $III''$ and $II_1$
\end{rem}

\begin{prop}\label{prop3}
Let us suppose that $\alpha$ is a morphism for $[ \cdot , \cdot ]$ and that 
 $[ \cdot , \cdot ]$ is a Lie algebra, then $[\cdot,
\cdot]_\alpha$ is a also a Lie algebra.
\end{prop}

\begin{proof} Let us consider again the equation
$$\text{ (\ref{lie_alpha}) } \mathfrak{J}_\alpha=\mathfrak{J}+
J^{I_1}_\alpha+J^{I_2}_\alpha+J^{I_3}_\alpha+J^{II}_\alpha+J^{II_2}_\alpha+J^{II_3}_\alpha.$$
It becomes, under the use of properties (\ref{property}) and
$S(\ref{property})$,
$$\mathfrak{J}_\alpha=\mathfrak{J}+J^{II}_\alpha-J^{II_1}_\alpha.$$
Now the hypothesis of being a morphism implies that
$J^{II}_\alpha=J^{II_1}_\alpha$, and hence the conclusion.
\end{proof}

Of course the preceding proposition does not involve Hom-algebras in its
formulation, but we think it is nevertheless useful to state it
because it admits the following generalization:

\begin{prop}\label{prop4}
Let us suppose that $([ \cdot , \cdot ], \alpha)$  is simultaneously
a $Hom^{II}$ and a $Hom^{II_1}$-algebra and that $[ \cdot , \cdot ]$ is  a Lie algebra, then $[\cdot,
\cdot]_\alpha$ is also a Lie algebra.
\end{prop}

\begin{proof} In the preceding proof, the hypothesis of being a morphism was used via the property
 $J^{II}_\alpha=J^{II_1}_\alpha$. But this last property is still satisfied under our
 new assumption $(0=0!)$, which ends the proof.
\end{proof}

 It is a good point to make a
little pause and have some side remarks. If one compares the two
previous statements, one sees that
 being simultaneously of Hom-algebra
types $\Phi$ and $\Phi'$ can be thought of as a replacement of the
notion morphism.
 But on the other hand the defect of being a morphism, i.e. the
 difference $C_\alpha (a,b):=[\alpha (a),\alpha(b)]-\alpha([a,b])$ is, in geometry, at the heart of the notion of curvature, and has many applications such as characteristic classes for example.
  This at once suggests to first introduce the defect of being
  simultaneously $Hom^\Phi$ and $Hom^{\Phi'}$:
$C^{\Phi}_\alpha (x,y,z):=\circlearrowleft_{x,y,z} [x,\alpha
([y,z])]-[x,[\alpha (y),\alpha (z)]]$ (of course one can similarly
define $C^{\star}_\alpha$ for the other types). A natural question
to investigate, would be to which extent the constructions in
homological algebra, based on the usual notion of curvature, can be
generalized in this framework. A second question would be to see if
these Hom-algebras carry any geometrical meaning, i.e. to search for
geometrical properties of tensors of type 1-1 on a manifold such
that the Lie algebra of vector fields, with this tensor, form a
Hom-algebra.

\vskip 1cm


\newpage

 \section{Hom-associative algebras.}
\label{homass}

In this part we start, applying the ideas of the previous part to
the associative
category, by getting in section \ref{section 2}
the list of different types of Hom-associative structures
one can consider. Subsequently, we focus our studies on \emph{unital} Hom-associative algebras since much more can be said about this special case than in general. In particular we state a hierarchy on such types of unital hom-algebras, i.e a complete classification of implications between them. We introduce hom-monoids which will give a usefull tool to build counterexamples.

We turn then in section \ref{proofs} to the proof of this hierarchy.
We start by proving implications that do hold and then give counter examples to the others.

\subsection{Types, unitality and hierarchy.}\label{section 2}

\subsubsection*{Types of Hom-associative algebras.}\label{astype}

Hom-associative algebras were introduced in \cite{MSI}, as examples
hom-structures. But there are more
 types of hom-associative algebras than the ones considered
 in \cite{MSI}. Analogously to hom-Lie algebras, one gets the following new types of hom-associative algebras:
\bigskip

\begin{tabular}{|c|c|c|c|}

\hline  &   & $II$  & $x\star \alpha (y\star z)=\alpha (x\star
y)\star z$
\\
\hline $I_1$  & $\alpha (x)\star (y\star z)=(x\star y)\star
\alpha(z)$ &$II_1$  & $x\star (\alpha (y)\star \alpha (z))=(\alpha
(x)\star \alpha (y))\star z$
\\
\hline $I_2$  & $x\star (\alpha (y)\star z)=(x\star \alpha (y))\star
z$ & $II_2$  & $\alpha (x)\star (y\star \alpha (z))=(\alpha
(x)\star y)\star \alpha (z)$\\
\hline $I_3$  & $x\star (y\star \alpha (z))=(\alpha (x)\star y)\star
z$ & $II_3$  & $\alpha (x)\star (\alpha (y)\star z))=(x\star \alpha
(y))\star
\alpha (z)$\\
 \hline
\end{tabular}




\begin{tabular}{|c|c|}

\hline    $III$  & $\alpha (x\star (y\star z))=\alpha ((x\star
y)\star z)$
\\
\hline    $III'$  & $\alpha (x)\star \alpha (y\star z)=\alpha
(x\star y)\star \alpha (z)$\\
 \hline $III''$  & $\alpha (x)\star
(\alpha (y)\star \alpha (z))=(\alpha (x)\star \alpha (y))\star
\alpha (z)$\\
 \hline
\end{tabular}
\bigskip

For precision, we give the following general definition:
\begin{df} \label{homring-definition}
Let $V$ be a set together with two binary operations $+:V \times V \rightarrow V$ and $\star:V \times V \rightarrow V$, one self-map $\alpha:V \rightarrow V$ and a special element $0 \in V$. Then $(V, +, \star, \alpha, 0)$ is called a \emph{hom-ring of type $T$} if 
\begin{itemize} 
\item $(V, +, 0)$ is an abelian group
\item the multiplication is distributive on both sides
\item $\alpha$ is an abelian group homomorphism
\item $\alpha$ and $\star$ satisfy the associativity condition corresponding to type $T$.
\end{itemize}
Hom-associative \emph{algebras} over a commutative ring $k$ are defined analogously by replacing any additivity conditions in the preceding definitions by corresponding conditions of $k$-linearity.
\end{df}
We remark that as in the associative case, a hom-ring may always be viewed as a hom-algebra with $\mathbb{Z}$ as base ring. Therefore, although we will in the sequel mainly talk about hom-algebras, of course only results that need conditions on properties of the base ring cannot be put in the hom-ring setting.\\
It seems that not much can be proven \emph{in general} about the relations among the various types of hom-associative algebras just introduced.
However, if additional conditions are imposed which restrict the range of algebras under consideration, the theory becomes much richer.
 One particularly natural condition which can be imposed is existence of a unit element.
We will see in section \ref{hierarch} that for \emph{unital} hom-associative algebras, the different types we defined can in some sense be partially ordered by increasing generality. 

Let us note that some of these types, namely $I_2, I_3$ and $II$ already appeared in \cite{FGII} but not under that name. One can reformulate these results (the meaning of unitality is precised below):

\begin{prop} \label{imassociative}Let $(V, \star, \alpha, 1)$ be a unital Hom-associative algebra of type $I_1$, then it is also of type $I_2$.
\end{prop}

\begin{df}
Let $(A, \star, \alpha)$ be a hom-associative algebra. Then $A$ is called \emph{left weakly unital} if $\alpha(x) = cx$ for some $c \in A$.
\end{df}
\begin{lem} \label{induced type I3 lemma}
Let  $(A, \star, \alpha, c)$ be a weakly left unital hom-associative algebra with weak left unit $c \in A$, bijective $\alpha$ and let $\beta := \alpha^{-1}$, then 
$(A, \star, \beta)$ is a hom-associative algebra of type $I_3$, it is also of type $II$.
\end{lem}

Proposition \ref{imassociative} is the transcription of proposition 1.1 of \cite{FGII}, while lemma \ref{induced type I3 lemma} concatenates lemmata 2.1 and 2.4 of \cite{FGII}.

\subsubsection*{Unitality.}\label{unit}

In the rest of the paper we assume $(V,\star, 1)$ to be unital. By
unitality we mean the usual notion in algebra, i.e the existence of
an element $1$ in V such that $1\star x=x\star 1$ for all $x$ in V.
The notion of an inverse of $x$ is also meaningful to some extent:
$x^{-1}$ is a left (resp. right) inverse of $x$ if it satisfies
$x^{-1}\star x=1$ (resp. $ x \star x^{-1}=1$). We call $x^{-1}$ an
inverse of $x$ if it is both a left inverse and a right inverse. There are known examples of non associative algebras with elements admitting different right and left inverses. Hence,  since \emph{any} bilinear map $\star: V \times V \rightarrow V$
induces a hom-associative structure on $V$ if we take $\alpha = 0$
as twisting homomorphism, there is no reason to assume that the
inverse of $x$ should be unique if defined. 

\subsubsection*{ Hierarchy on types of unital Hom associative
algebras.}\label{hierarch}


 The hierarchy alluded to in the introduction can for the relations among first and second types be summarized as follows:

\begin{prop}\label{resume}
Let $(V,\star, 1)$ be unital, then one has the following relations
between the types of $(V,\star, 1, \alpha)$:\\ a) $II_1\Leftarrow
II\Leftrightarrow I_1 \Rightarrow I_3 \Rightarrow \{I_2,II_2,II_3 \
and \ II_1 \}$\\
b) $I_2\Rightarrow \{II_1\Leftrightarrow II_3\}.$
\end{prop}

There are no relations among \{types III\} and from \{types III\} to
\{types I and II\}. The relations from \{types I and II\} to \{types
III\} are given by:

\begin{prop}\label{resumeII}
Let $(V,\star, 1)$ be unital, then one has the following relations
between the types of $(V,\star, 1, \alpha)$:\\ a) $I_1\Leftarrow
III, III', III''$\\
b) $I_3\Rightarrow
III''$\\
c) $I_2\Rightarrow
III''$\\
d) $II_2\Rightarrow
III''$\\
e) $II_1, II_3\Rightarrow
III''$\\
\end{prop}

 The equivalence
of types $I_1$ and $II$ will been given in Proposition
\ref{IeqII}. Proposition \ref{basic} gives that type
$I_1$ implies type $I_3$ and that type $II$ implies type $II_1$. Proposition
\ref{I_3impl} states that $I_3$ implies $I_2,II_2,II_3$ and $II_1$.
We prove in Proposition \ref{exot} the exotic implication of
Proposition \ref{resume} b).

Proposition \ref{resumeII} is the conjunction of Propositions
\ref{prop6}, \ref{prop7} and \ref{prop5}.

\subsubsection*{Hom-monoids}\label{counterexamples} We will now introduce \emph{hom-monoids} and give a short
discussion of their relation to hom-algebras. The main motivation is
that hom-monoids will be useful in the construction of
counterexamples to relations between the different types of
hom-algebras.
\begin{df}
A \emph{hom-monoid} of type $I$ is a set $S$ together with a binary operation $\star:S \times S \rightarrow S$, a special element $1 \in S$ and a map $\alpha:S \rightarrow S$ such that the following axioms are fulfilled:
\begin{eqnarray*}
1 \star x & = & x \star 1 = x\\
\alpha(x) \star (y \star z) & = & (x \star y) \star \alpha(z).
\end{eqnarray*}
Similarly, we introduce for each type of hom-associative algebra defined previously the corresponding type of hom-monoid. If we do not specify a type, Type $I$ will by default be implied.
\end{df}
The following example is clear:
\begin{ex} Let $(V, \star, +, \alpha, 1)$ be a hom-algebra of type $T$. Then the multiplicative structure $(V, \star, \alpha, 1)$ is a hom-monoid also of type $T$.
\end{ex}
In the other direction, one has the following remark:
\begin{rem} Let $k$ be a commutative ring and let $(S, \tilde{\star}, \tilde{\alpha}, 1)$ be a hom-monoid of type $T$. Let then $V$ be the free $k$-module over $S$ and define $\alpha:V \rightarrow V$ and $\star:V \times V \rightarrow V$ by linear extension of $\tilde{\alpha}:S \rightarrow S$ respectively $\tilde{\star}:S \times S \rightarrow S$ to $V$. Then $(V, \star, \alpha, 1)$ is a unital hom-associative algebra of type $T$. We denote the hom-algebra so constructed from a hom-monoid $S$ by $k[S]$.
\end{rem}
\begin{proof} By construction, $\alpha$ is linear and $\star$ is bilinear. Using the distributive laws, one verifies easily that type $T$ hom-associativity of $(V, \star, \alpha)$ follows from the corresponding property on the generating set $S$. Unitality of $V$ is clear. This concludes the proof.
\end{proof}
With respect to exploring the relations between different types of hom-associative algebras, the preceding remark and example show that if type $T_1$ subsumes type $T_2$ in the context of unital hom-algebras, then the same holds in the context of hom-monoids and  vice versa.\\
To obtain from hom-monoids examples of hom-algebras which can be written down in a particularly concise way, it will be useful to also set the following definition:
\begin{df}
A hom-monoid $(S, \star, \alpha, 1)$ of type $T$ is called a hom-monoid of type $T$ \emph{with zero} if there is an element $0 \in S$, $0 \neq 1$, such that $0 \star x = x \star 0 = 0$ for all $x \in S$ and $\alpha(0) = 0$.
\end{df}
Since as with unital hom-algebras of the original type $I_1$, the twisting map $\alpha$ is also in a hom-monoid of type $I_1$ automatically multiplication with some element inside the structure, one could drop the condition $\alpha(0) = 0$ in this case from the definition of a hom-monoid with zero. However, for most of the other types it is necessary to impose it separately.\\
In any case, if $V$ is a hom-algebra over some commutative ring $k$ of type $T$ constructed from a hom-monoid $S$ with zero of the same type, and if $0_S$ denotes the zero element in $S$, then the submodule $I := k \cdot 0_S$ becomes a hom-ideal of $k[S]$. One can kill this submodule by passing to the appropriate factor algebra $A := k[S]/I$. We note that $A$ still contains a copy of $S$. In particular, if $S$ was (not) of type $T$, then $A$ will be (not) of type $T$.\\

\subsection{Proof of the hierarchy} \label{proofs}

The rest of this paper
is devoted to the exploration of the logical relationships between the different types
of unital hom-algebras. We start by establishing relations that do hold. Counter examples to the other relations are given at the end of this part. 
\subsubsection*{Equivalence between types $I_1$ and $II$}

We start by a lemma which contains the main basic properties
allowing computations with Types $I_1$ and $II$. This lemma was proved in \cite{FGII} (lemma 1.1) under the assumption of being of type $I_1$. We recall it and extend it under the assumption of being of type $II$.

\begin{lem}\label{lemma1}
Let $(V,\star,\alpha,1)$ be of type $I_1$ or $II$. One has for all
$x,y$ in $V$:
\begin{enumerate}
\item[\"{a})] $\alpha(x)\star y=x\star \alpha (y),$
\item[\"{e})] $x\star \alpha (1)=\alpha (x),$
\item[\"{\i})] $\alpha (x\star y)=x\star \alpha (y).$
\end{enumerate}
\end{lem}

\begin{proof} The
proof of the two first points, assuming the Hom-algebra to be of
type $II$ works along exactly the same lines as the proof of lemma 1.1 in \cite{FGII} and is left to the
reader. To prove \"{\i}), simply apply the definition of Hom algebra
of type $II$ to the
triple $(x,y,1)$: $\alpha(x\star y)\star 1\overset{II}{=}x\star \alpha(y\star 1).$
\end{proof}

 The desired equivalence is then a simple corollary:

\begin{prop}\label{IeqII}
Hom-associative algebras of types $I_1$ and $II$ are equivalent.
\end{prop}

\begin{proof} The contemplation of the following square gives the proof.
$$
  \begin{array}{ccc}
    x\star\alpha(y\star z) & \overset{II}{=} & \alpha(x\star y)\star z \\
    \text{\"{a)}}\shortparallel &   & \text{\"{a)}}\shortparallel \\
     \alpha(x)\star(y\star z) & \overset{I_1}{=} & (x\star y)\star \alpha(z).  \\
  \end{array}$$
\end{proof}

\subsubsection*{Type implications from $I_1$}  Now, we
concern ourselves with the types subsumed by the equivalent types
$I_1$ and $II$.
\begin{prop}\label{basic}
Let $(V, \star, \alpha, 1)$ be a unital Hom-associative algebra of
type $I_1$, then \\
a) it is also of type $I_3$;\\
b) it is also of type $II_1.$
\end{prop}
The proof of proposition \ref{basic} requires the following equalities taken from proposition 1.1 of \cite{FGII}.

\begin{lem}\label{imassoc} Let $(V, \star, \alpha, 1)$ be a unital Hom-associative algebra of type
$I_1$. One has for all $x,y$ and $z$ in $V$:
\begin{enumerate}
\item[1)] $\alpha(x)\star (y\star z)=(\alpha(x)\star y)\star z,$
\item[2)]  $x\star (y\star \alpha(z))=(x\star y)\star \alpha(z)$
\end{enumerate}
\end{lem}

\begin{proof}
The proof of a) comes from contemplation of the following square:
$$
  \begin{array}{ccc}
    \alpha(x)\star(y\star z) & \overset{I_1}{=} & (x\star y)\star \alpha(z) \\
    \text{\ref{imassoc} 1)}\shortparallel &   & \text{\ref{imassoc} 2)}\shortparallel \\
    ( \alpha(x)\star y)\star z & \overset{I_3}{=} & x\star (y\star \alpha(z)).  \\
  \end{array}$$\\

There exists a direct proof of b), but it can also be seen from the
chain $I_1 \overset{a)}{\Rightarrow} I_3
\overset{\ref{I_3impl}}{\Rightarrow}
\{I_2,II_3\}\overset{\ref{exot}}{\Rightarrow} II_1$ which will be
proven in the following two sections.
\end{proof}

 \vskip 1cm

\subsubsection*{Type implications from $I_3$.}\label{I3}

The main result of this section is

\begin{prop}\label{I_3impl}
Let $(V, \star, \alpha, 1)$ be a unital Hom-associative algebra of
type $I_3$, then \\
a) it is also of type $I_2$;\\
b) it is also of type $II_2$;\\
c) it is also of type $II_3.$
\end{prop}

Its proof is based on the following two basic properties:

\begin{lem}\label{bis}
Let $(V, \star, \alpha, 1)$ be a unital Hom-associative algebra of
type $I_3$, then $\forall x, y \in V,$\\
$\alpha)$ $\alpha(x)=x\star\alpha(1)$;\\
$\beta)$ $x\star\alpha(y)=\alpha(x)\star y.$
\end{lem}

\begin{proof}[Proof of the lemma]
 One proves $\alpha)$ by applying the definition of
$Hom^{I_3}$-associativity to the triple $(x,1,1)$:
$$x\star\alpha(1)=x\star(1\star \alpha(1))\overset{I_3}{=}(\alpha(x)\star1)\star 1=\alpha(x).$$

The proof of $\beta)$ is obtained by considering the triple
$(x,y,1)$ by the use of a):
$$x\star\alpha(y)\overset{a)}{=}x\star(y\star\alpha(1))\overset{I_3}{=}(\alpha(x)\star y)\star1)=\alpha(x)\star
y.$$\\
\end{proof}

One should resist the temptation to deduce lemma \ref{lemma1} from proposition \ref{basic} and lemma \ref{bis} since proposition \ref{basic} relies itself on lemma \ref{lemma1}.
We now turn on to the proof of Proposition \ref{I_3impl}:

\begin{proof} We show each statement in turn:\\
a) $$
  \begin{array}{ccc}
    x\star(\alpha(y)\star z) & \overset{I_2}{=} & (x\star \alpha(y))\star z \\
    \beta)\shortparallel &   & \shortparallel \beta)\\
     x\star (y \star \alpha(z)) & \overset{I_3}{=} & (\alpha(x)\star y)\star z).  \\
  \end{array}$$\\
b)
$$
  \begin{array}{ccc}
    (\alpha(x)\star y)\star \alpha(z) & \overset{II_2}{=} & \alpha(x)\star (y\star \alpha(z)) \\
I_3\shortparallel &   & \shortparallel I_3\\
     x\star (y \star \alpha(\alpha(z))) &  & (\alpha(\alpha(x))\star y))\star z.  \\
    \beta)\shortparallel &   & \shortparallel \beta)\\
     x\star (\alpha(y) \star \alpha(z)) & \overset{I_3}{=} & (\alpha(x)\star \alpha(y))\star z.  \\
  \end{array}$$\\
c)
$$
  \begin{array}{ccc}
    \alpha(x)\star (\alpha(y)\star z) & \overset{II_3}{=} & (x\star \alpha(y))\star \alpha(z) \\
\beta)\shortparallel &   & \shortparallel \beta)\\
     \alpha(x)\star (y \star \alpha(z)) &  & (\alpha(x)\star y)\star \alpha(z).  \\
I_3\shortparallel &   & \shortparallel I_3\\
     (\alpha(\alpha(x))\star y) \star z &  &  x\star (y \star \alpha(\alpha(z))).  \\
    \beta)\shortparallel &   & \shortparallel \beta)\\
     (\alpha(x)\star \alpha(y)) \star z & \overset{I_3}{=} & x\star(\alpha(y)\star \alpha(z)).  \\
  \end{array}$$\\
\end{proof}


\subsubsection*{Exotic implications.}\label{exotic}

We call these implications exotic since, contrary to the previous
ones, they involve two types of hom-algebras in the
assumptions.

\begin{prop}\label{exot}
Let $(V, \star, \alpha, 1)$ be a unital Hom-associative algebra of
type $I_2$, then it is of type
$II_3$ if and only if it is of type $II_1.$
\end{prop}
\begin{proof} The following calculation shows our claim:
$$
  \begin{array}{ccc}
    x\star(\alpha(y)\star \alpha(z)) & \overset{II_1}{=} & (\alpha(x)\star \alpha(y))\star z \\
    I_2\shortparallel &   & \shortparallel I_2\\
    ( x\star \alpha(y))\star \alpha(z) & \overset{II_3}{=} & \alpha(x)\star (\alpha(y)\star z).  \\
  \end{array}$$\\
\end{proof}









\subsubsection*{Implications from types of families I and II to types of family III.}

There are not many relations between these types, except for $III"$
which is weaker than almost all the other types and $I_1$ which is
stronger than all the other types.

\begin{prop}\label{prop6} One has the following implications\\
a) a hom-algebra of type $I_1$ is necessarily of type
$III$,\\
b)  a hom-algebra of type $I_1$ is necessarily of type
$III'$,\\
c)  a hom-algebra of type $I_1$ is necessarily of type
$III''$.
\end{prop}

\begin{proof}
a) is Proposition \ref{imassociative}, $(4)$.\\
b) $$ \begin{array}{ccc}
    \alpha(x)\star\alpha(y\star z) & = &  \alpha(x\star y)\star \alpha(z) \\
   \text{\it lemma }\ref{lemma1}\shortparallel  &   & \shortparallel \text{\it lemma }\ref{lemma1}\\
  \alpha(x)\star (\alpha(y)\star z)   & \overset{I_1}{=} & (x\star \alpha(y))\star \alpha(z).  \\
  \end{array}$$\\
{\it Hom III'}-associativity in c) is obtained applying proposition
\ref{imassociative} 1) to the triple $\{x, \alpha(y),\alpha(z)\}$.\\
\end{proof}

The last point of the previous proposition can be refined (and
implied) by:

\begin{prop}\label{prop7} A hom-algebra of type $I_3$ is necessarily of type
$III''.$
\end{prop}

\begin{proof}
By Proposition \ref{I_3impl} a), $I_3\Rightarrow I_2$. But by
Proposition \ref{prop5} a) below, $I_2\Rightarrow III'$.
\end{proof}

There are no other implications in this direction from $I_3$, as
shown by the following counterexample:

$I_3 \nRightarrow III, III':$ $(e_2\cdot e_2):=e_1, (e_3\cdot e_3):=e_3$;
$\alpha(e_1)=\alpha(e_3)=e_3$.\\

In particular, since $I_3$ implies $\{II_1, II_2, II_3, I_2\}$, the
previous counter examples are also counter examples to $\{II_1,
II_2, II_3, I_2\}\Rightarrow III$ and $\{II_1, II_2, II_3,
I_2\}\Rightarrow III'$. But can $III''$ be implied by one of these
types?

\begin{prop}\label{prop5}
One has the following implications:

a) a hom-algebra of type $I_2$ is necessarily of type $III''$\

b) a hom-algebra of type $II_2$ is necessarily of type $III''$\

c) a hom-algebra of types $II_1$ and $II_3$ is necessarily of type $III''$.\

\end{prop}

 The hard point to prove is c), and we will need for it the
following

\begin{lem}\label{lemma3} For a hom-algebra $V$  which is of types $II_1$ and $II_3$, one has
$\forall x, y \in V$ :

a) $\alpha (1)\star \alpha (x)= \alpha(x)\star\alpha (1) $\

 b) $ \alpha (\alpha(x))\star\alpha (1)=
\alpha(x)\star\alpha (\alpha (1))$\

c)
$\alpha(x)\star(\alpha(1)\star\alpha(y))=(\alpha(x)\star\alpha(1))\star\alpha(y)$\

 d) $ \alpha
(\alpha(x))\star\alpha (y)= \alpha(x)\star\alpha (\alpha (y)).$
\end{lem}

\begin{proof}[Proof of Lemma \ref{lemma3}] We start with (a):
$$ \begin{array}{ccc} \alpha (1)\star \alpha (x) & \overset{}{=}
& \alpha(x)\star\alpha (1)\\
\shortparallel  &   & \shortparallel  \\
 (\alpha (1)\star \alpha (x))\star 1 & \overset{II_1}{=}
& 1\star (\alpha(x)\star\alpha (1))\\
\end{array}$$

b) $$ \begin{array}{ccc}
 \alpha(\alpha(x))\star\alpha (1) & \overset{}{=}
&1\star(\alpha (x)\star \alpha(\alpha (1)))  \\
\shortparallel  &   & \shortparallel II_1 \\
(\alpha(\alpha(x))\star \alpha (1))\star 1 &
&(\alpha(1)\star \alpha(x))\star \alpha (1)   \\
II_1\shortparallel  &   & \shortparallel a) \\
\alpha(x)\star (\alpha(1)\star \alpha (1)) & \overset{II_3}{=}
&(\alpha(x)\star \alpha(1))\star \alpha (1)  \\

\end{array}$$

 c) $$ \begin{array}{ccc}
\alpha(x)\star(\alpha(1)\star\alpha(y)) & \overset{}{=} & (\alpha(x)\star\alpha(1))\star\alpha(y)\\
II_1\shortparallel  &   & \shortparallel  II_1 \\
(\alpha(\alpha(x))\star \alpha(1))\star y & & x\star (\alpha(1)\star (\alpha(\alpha(y)))) \\
b)\shortparallel  &   & \shortparallel  b) \\
(\alpha(x)\star \alpha(\alpha(1)))\star y & \overset{II_1}{=} &
x\star (\alpha(\alpha(1))\star (\alpha(y)))
\end{array}$$

d)$$ \begin{array}{ccc}
\alpha(\alpha(x))\star\alpha(y) & = & \alpha(x)\star\alpha(\alpha(y))\\
\shortparallel  &   & \shortparallel  \\
(\alpha(\alpha(x))\star\alpha(y))\star 1 &  & 1 \star(\alpha(x)\star\alpha(\alpha(y)))\\
II_1\shortparallel  &   & \shortparallel  II_1 \\
\alpha(x)\star(\alpha(y)\star \alpha(1)) &  & (\alpha(1) \star\alpha(x))\star\alpha(y)\\
a)\shortparallel  &   & \shortparallel a) \\
\alpha(x)\star(\alpha(1)\star\alpha(y)) & \overset{c)}{=} &
(\alpha(x)\star\alpha(1))\star\alpha(y)
\end{array}$$
\end{proof}

\begin{proof}[Proof of Proposition \ref{prop5}]

One proves a) by applying the definition of type $I_2$ to the triple
$\{ \alpha(x), y, \alpha(z)\}$, and one proves b) by applying the
definition of type $II_2$ to the triple $\{ x, \alpha(y), z\}$. The
proof of c) comes from contemplation of the following diagram:

$$ \begin{array}{ccc}
     (\alpha(x)\star\alpha(y))\star \alpha(z)& \overset{III''}{=} & \alpha(x)\star (\alpha(y)\star\alpha(z)) \\
II_1\shortparallel  &   & \shortparallel  II_1 \\
 x\star (\alpha(y)\star \alpha(\alpha(z))) & & (\alpha(x)\star\alpha(y))\star
 \alpha(z)\\
   \text{\it lem }\ref{lemma3}\shortparallel  &   & \shortparallel \text{\it lem }\ref{lemma3}\\
  x\star (\alpha(\alpha(y))\star \alpha(z)) &  \overset{II_3}{=} &  (\alpha(x)\star \alpha(\alpha(y)))\star z .  \\
  \end{array}$$

\end{proof}

\subsubsection*{Counterexamples to intertype relations.}

We give now a list of counterexamples to inter-type relations which do \emph{not} hold. All of these counterexamples are constructed by use of the technique of hom-monoids discussed in section \ref{section 2} above. For each example, we first describe what the example is supposed to show. Here, $T_1, T_2, \ldots, T_n \not\Rightarrow T$ means that a hom-monoid which is of simultaneously of types $Type_i$ is not necessarily also type $T$. This is followed by relations between elements of a counterexample hom-monoid of the requisite types. The elements of the hom-monoid structures in question are denoted $e_1, \ldots, e_n$, where $e_1$ is supposed the unit element. All hom-monoid structures are with zero, but the zero element is outside the set $\{e_1, \ldots, e_n\}$. When we do not give a product $e_i e_j$, this means that $e_i e_j = 0$, except when $i = 1$ or $j = 1$, in which case the product is prescribed by the requirement that $e_1$ be the unit of the hom-monoid. The hom-monoids in question are not supposed to have elements $e_i$ with higher index than those appearing in the relations we supply. Likewise, we give only nonzero values of $\alpha$. Note that the first three examples use the fact that not every associative structure is also hom-associative of all types.\\
\begin{enumerate}

\item $I_2 \nRightarrow I_3$:
$\alpha(e_2)=e_1$.


\item $I_2 II_2\nRightarrow II_3$:
  $\alpha(e_1)=e_1$.

\item $II_1II_2 II_3I_2\nRightarrow I_3$:
$\alpha(e_2)=e_2$.



\item $I_3 \nRightarrow I_1$:
$(e_2\cdot e_2):=e_1$; $\alpha(e_2)=\alpha(e_3)=e_3$.

\item $I_2 \nRightarrow II_2$:
$(e_2\cdot e_2):=e_2, (e_2\cdot e_3):=e_2$; $\alpha(e_1)=\alpha(e_2)=e_3$.

\item $II_2 \nRightarrow I_2$:
$(e_2\cdot e_3)=e_3$; $\alpha(e_1)=e_2$.

\item $II_1II_2 \nRightarrow II_3$:
$(e_2\cdot e_3)=(e_3\cdot e_2)=e_1$; $\alpha(e_1)=e_3$.

\item $II_1\nRightarrow II_2$:
$(e_2\cdot e_3)=e_1$; $\alpha(e_1)=e_2$.

\item $II_1II_2 II_3\nRightarrow I_2$:
$(e_2\cdot e_3)=e_1, (e_3 \cdot e_2)=e_1$; $\alpha(e_2)=e_3$.

\item $II_2II_3\nRightarrow II_1$:
$(e_2\cdot e_2)=e_1, (e_3\cdot e_3)=e_2$;
$\alpha(e_1)=e_3$.



\item $I_2II_1 II_3\nRightarrow II_2$:
$(e_3\cdot e_2)=e_4, (e_4\cdot e_3)=e_2$; $\alpha(e_1)=e_3$. The next three examples show the absence of relations between Hom types $III, III'$ and $III''$:



\item $III, III'' \nRightarrow III'$:
$\alpha(e_2)=e_1$.


\item $III, III' \nRightarrow III''$:
$(e_2\cdot e_2):=e_3,(e_3\cdot e_2):=e_2$;$\alpha(e_1)=e_2$.

\item $III', III'' \nRightarrow III$:
$(e_2\cdot e_2):=e_3, (e_3\cdot e_2):=e_3$; $\alpha(e_3)=e_3$.\\
Lastly, we give an example that shows that the order three types do not subsume any of the other types even when taken together:
\item $III, III', III'' \nRightarrow I_2, II_1, II_2, II_3$: $e_2 \cdot e_2 = e_1, e_2 \cdot e_3 = e_1, e_3 \cdot e_2 = e_2, e_3 \cdot e_3 = e_1$ and $\alpha(e_1) = \alpha(e_2) = \alpha(e_3) = e_3.$
\end{enumerate}



\vskip 1cm

 \noindent {\bf Acknowledgements.}  The authors were supported in their work by research grant R1F105L15 of the University of Luxembourg (Martin Schlichenmaier).\\
This work was aided by use of the first order logic theorem prover Prover9 and the finite countermodel searcher Mace4. In particular, Mace4 was very useful in the practical application of the Hom-monoid technique for counterexample building in the proof of the ``negative'' parts of the Hom-associative type hierarchy. Also some of the examples in the Lie case are due to Mace4. Prover9 was useful for quick verification of many parts of the associative hierarchy, although we mostly did not use the proofs offered.\\
The first author wants to thank the Max Planck Institute for Mathematics in Bonn for excellent working conditions.\\
We also thank Abdenacer Makhlouf and Donald Yau for usefull comments on a preliminary draft.
This work has been published thanks to the financial support of the "Fonds National de la Recherche" of Luxembourg.
\vskip 1cm


\end{document}